\documentclass[12pt]{amsart}
\usepackage{amsmath,amsfonts,amssymb,amsthm,epsfig}
\usepackage{latexsym}
\usepackage{amssymb}
\usepackage{graphicx}
\usepackage{amssymb}

\usepackage{hyperref}

\usepackage[margin=1.0in]{geometry}

\def\({\left(}
\def\){\right)}
\def\[{\left[}
\def\]{\right]}
\def\<{\left\langle}
\def\>{\right\rangle}

\newtheorem{theorem}{Theorem}[section]

\newtheorem{corollary}[theorem]{Corollary}

\newtheorem*{question*}{Question}
\newtheorem*{conjecture*}{Conjecture}

\DeclareMathOperator{\supp}{supp}

\begin{document}

\author[A. Bufetov]{Alexander Bufetov}
\email{bufetov@mi.ras.ru}
\address{Laboratoire d'Analyse, Topologie, Probabilit\'es, CNRS, Marseille}
\address{The Steklov Institute of Mathematics, Moscow}
\address {The Institute for Information Transmission Problems, Moscow}
\address{National Research University Higher School of Economics, Moscow}
\address{The Independent University of Moscow}
\address{Rice University, Houston, USA}

\author[S. Mkrtchyan]{Sevak Mkrtchyan}
\email{sevakm@math.cmu.edu}
\address{Carnegie Mellon University, Pittsburgh, USA}

\author[M. Shcherbina]{Maria Shcherbina}
\email{shcherbi@ilt.kharkov.ua}
\address{Institute for Low Temperature Physics Ukr. Ac. Sci., Kharkov, Ukraine}

\author[A. Soshnikov]{Alexander Soshnikov}
\email{soshniko@math.ucdavis.edu}
\address{University of California at Davis, Davis, USA}

\title[Entropy for beta random matrix ensembles]{Entropy and the Shannon-McMillan-Breiman theorem for beta random matrix ensembles}

\begin{abstract}
We show that beta ensembles in Random Matrix Theory with generic real analytic potential have the asymptotic equipartition property.
In addition, we prove a Central Limit Theorem for the density of the eigenvalues of these ensembles.
\end{abstract}

\maketitle

\section{Introduction}

In this paper we study asymptotic properties of the density functions of certain measures known as beta ensembles
that arise in Random Matrix Theory.  Namely, we consider probability distributions in $\mathbb{R}^N$ of the form

\begin{align}
\label{meas:Hermitian}
P^{\beta}_N(\lambda_1,\lambda_2,\dots,\lambda_N)&=
\frac1{Z_N(\beta)}\prod_{i<j}|\lambda_i-\lambda_j|^\beta\prod_{i=1}^Ne^{-\beta N \* V(\lambda_i)/2},
\end{align}
where the potential $V$ is a real analytic function satisfying the growth condition
\begin{equation}
\label{potential}
V(\lambda) \geq 2\*(1+\varepsilon)\*\log(1+|\lambda|)
\end{equation}
for all sufficiently large $\lambda.$

Beta ensembles have attracted significant interest in recent years (see e.g. \cite{BouErdYau1}, \cite{BouErdYau2}, 
\cite{DumEdeBetaModel}, \cite{DumEdeGlobal}, \cite{DumPaq}, \cite{ErdYau},
\cite{Pop},  \cite{RamRidVir}, \cite{S:11}, \cite{S:12},
\cite{SosWon}, \cite{ValVir}, \cite{WongLocal}, and references therein).
Below, we briefly mention four classical beta ensembles, namely the Hermite (Gaussian) beta ensemble, the Circular beta ensemble, the 
Laguerre (Wishart) beta ensemble, and the Jacobi (MANOVA) beta ensemble.

\begin{itemize}
\item\textit{Hermite beta ensemble: }
Let $\lambda_1,\dots,\lambda_N\in\mathbb{R}$ be random variables with joint density function with respect to the Lebesgue measure given by
\begin{align}
\label{meas:beta}
P^{Her,\beta}_N(\lambda_1,\lambda_2,\dots,\lambda_N)&=
\frac1{Z^{Her}_N(\beta)}\prod_{i<j}|\lambda_i-\lambda_j|^\beta\prod_{i=1}^Ne^{-\frac{\beta N}{4}\lambda_i^2}.
\end{align}
For $\beta=1, 2,$ and $4,$ the distribution (\ref{meas:beta}) is known as
the joint distribution of the eigenvalues of a random matrix from the Gaussian
Orthogonal Ensemble (GOE), Gaussian Unitary Ensemble (GUE), and Gaussian Symplectic Ensemble (GSE) respectively.  Dumitriu and Edelman
(see \cite{DumEdeBetaModel}, \cite{AndGuiZei}) introduced tridiagonal real symmetric
random matrices with i.i.d. centered Gaussian random variables on the diagonal and $\chi$ distributed independent random variables on
the upper and lower sub-diagonals such that the joint distribution of the eigenvalues is given by (\ref{meas:beta}) for arbitrary $\beta>0.$

The next three classical ensembles do not formally belong to the class (\ref{meas:Hermitian}) since the particles are distributed, respectively, on
the unit circle, positive half-line, and the interval $[-1,1].$

\item\textit{Circular beta ensemble: }
Let $\lambda_1,\dots,\lambda_N\in[0,2\pi]$ be random variables with the joint density function with respect to the Lebesgue measure given by
\begin{align}
\label{meas:circbeta}
P^{Cir,\beta}_N(\lambda_1,\lambda_2,\dots,\lambda_N)&=\frac{1}{Z^{Cir}_N(\beta)}\prod_{k<j}|e^{i\lambda_k}-e^{i\lambda_j}|^\beta.
\end{align}

\item\textit{Laguerre ensemble: }
Let $\lambda_1,\dots,\lambda_N\in[0,\infty)$ be random variables with the joint density function with respect to the Lebesgue measure given by
\begin{align}
\label{meas:lag}
P^{Lag,\beta}_N(\lambda_1,\lambda_2,\dots,\lambda_N)&=\frac{1}{Z^{Lag}_N(\beta)}\prod_{i<j}|\lambda_i-\lambda_j|^\beta
\prod_{j=1}^N\lambda_j^{\alpha-1}e^{-\beta N \lambda_j},
\end{align}
where $\alpha>0$.

\item\textit{Jacobi ensemble: }
Let $\lambda_1,\dots,\lambda_N\in[-1,1]$ be random variables with the joint density function with respect to the Lebesgue measure given by
\begin{align}
\label{meas:jac}
P^{Jac,\beta}_N(\lambda_1,\lambda_2,\dots,\lambda_N)&=\frac{1}{Z^{Jac}_N(\beta)}
\prod_{i<j}|\lambda_i-\lambda_j|^\beta\prod_{j=1}^N(1-\lambda_j)^{\mu-1}(1+\lambda_j)^{\nu-1},
\end{align}
where $\mu,\nu>0$.

\end{itemize}

As in the Hermite case, the values $\beta=1,2,$ and $4$ (with $\alpha=\frac{\beta}{2}\*(n-m+1)-1$ in the Laguerre case and
$\mu= \frac{\beta}{2}\*(n_1-m+1)-1, \ \mu= \frac{\beta}{2}\*(n_2-m+1)-1$ in the Jacobi case),
correspond to classical ensembles of random matrices (see e.g. \cite{mehta2004random}, \cite{DumEdeBetaModel}). In the Laguerre case, the random
matrix ensemble for arbitrary positive $\beta$ was constructed in  \cite{DumEdeBetaModel}.
For Jacobi and Circular beta ensembles, three- and five-diagonal matrix models were derived by Killip and Nenciu in \cite{KilNen04}.

It is well known (see \cite{BG:11}, \cite{JohDuke1998}, \cite{Deift})
that if the potential $V$ in (\ref{meas:Hermitian})
is sufficiently smooth (e.g. when $V'$ is H\"older continuous), then there exists an equilibrium measure $\mu^V$
which is absolutely continuous with
respect to the Lebesgue measure and has compact support. Let us denote its density by $\rho^V.$ The equilibrium measure maximizes the functional
\begin{equation}
\label{energy}
\mathcal{E}_V(\mu):=\int \int \log|x-y|\* d\mu(x)\*d\mu(y) -  \int V(x) \* d\mu(x)
\end{equation}
over the space of the probability measures on  $\mathbb{R}.$  Note that $\mathcal{E}_V(\mu),$ up to a factor of $-1,$ coincides with the energy
functional.
We denote by $\mathcal{E}[V]:=\mathcal{E}_V(\mu^V)$ the maximum value of the functional attained at the equilibrium measure  $\mu^V.$
The marginal density of the ensemble (\ref{meas:Hermitian}),
$$\rho_{1, \beta}^{(N)}(\lambda)=\int_{\mathbb{R}^{N-1}} P^{\beta}_N(\lambda,\lambda_2,\dots,\lambda_N) \*d\lambda_2\ldots d\lambda_N,$$
weakly converges to the equilibrium density $\rho^V.$  In addition, the support of $\mu^V$ and the density $\rho^V$ are uniquely determined by the
Euler-Lagrange variational equations

\begin{align}
\label{EL1}
& 2\* \int \log|x-y| \*d\mu^V(y) -V(x) \leq l, \ \ x\in \mathbb{R}, \\
\label{EL2}
& 2\* \int \log|x-y| \*d\mu^V(y) -V(x) =l, \ \ x: \rho^V(x)>0.
\end{align}
In the Hermite case $(V(x)=x^2/2),$ the equilibrium measure is the Wigner semicircular distribution
$\mu_{sc}(dx)=\frac{1}{2\pi}\*\sqrt{4-x^2} \* 1_{[-2,2]}(x)\*dx.$

In addition to (\ref{potential}), we
assume throughout this paper that the potential $V$ is real analytic and generic 
(see Theorem \ref{thm:AEP} for details).  In particular, the support of the equilibrium measure consists of finitely many intervals (see e.g. 
\cite{APS}).
We prove that the
logarithm of the density, after appropriate normalization, converges to a constant almost surely.
This is called the asymptotic equipartition property.
The precise statement is given in Theorem \ref{thm:AEP}.
By analogy with the Shannon--McMillan-Breiman theorem, the limiting constant can be interpreted
as the entropy of the corresponding measure. A corollary of this result states that in the limit when the size of the considered matrices goes to
infinity, there is a set of measure almost $1$ such that at all the points of this set the density is almost the same.

Analogous results for the
Plancherel measure were conjectured to be true by Vershik and Kerov \cite{VK85} and proved in \cite{Bu}. For a one parameter
deformation of the Plancherel measure, called the Schur--Weyl measures, equivalent results were obtained in \cite{M2,M3EntSW}. Note, that in the
case
of the Plancherel and Schur-Weyl measures it is only known that the normalized logarithm of the density converges to a constant in probability. The
question of almost sure convergence for those measures is still open.

In addition to the asymptotic equipartition property we prove the Central Limit Theorem for the logarithm of the density (see Theorem \ref{thm:CLT}).
For general random matrix models Borodin and Serfaty have studied asymptotic properties of a similar statistic, which they call ``renormalized
energy''
\cite{BorSer2011}. They explicitly calculate the limit of the expectation of the ``renormalized energy'' in the case of $\beta$-sine processes for
$\beta=1,2,4$ and some $2$-dimensional point processes. The result of Theorem \ref{thm:CLT} in the special
case of the circular beta ensemble appears in
\cite{BorSer2011}. The Central Limit Theorem for random variables of the form $\sum_i g(\lambda_i)$, where $\lambda_i$ are the eigenvalues and
$g$ is a bounded continuous function, was obtained by Johansson for the circular ensemble \cite{JohSzego1988} and later for ensembles with more
general
potentials \cite{JohDuke1998} (see also \cite{Pastur} for results on non-Gaussian limiting fluctuation of linear statistics in the multi-cut case).
A large deviation principle for spectral measures of certain classes of beta ensembles of random matrices
was established in \cite{BenArGui} (see also recent results by Forrester for the classical beta ensembles in \cite{Fo:12a} and \cite{Fo:12b}).

For the classical beta ensembles (Hermite, Circular, Laguerre, and Jacobi beta ensembles) our results follow from
the Selberg Integral \cite{Selberg} since it allows us to  obtain explicit formulas for the Laplace transform of the normalized
logarithm of our density.  We obtain the pointwise convergence and the Central Limit Theorem from the asymptotics of the exponential moment.
In the general case, we use recent results \cite{S:11}, \cite{BG:11}, \cite{S:12} about the asymptotic expansion of the partition function of
(\ref{meas:Hermitian}). We note that 
the use of the Selberg integral for deriving fluctuation formulas for linear statistics in the classical beta
ensembles was first done by Baker and Forrester in \cite{BakerFor}.
\subsection{Acknowledgements}
We are grateful to Alexei Borodin and Kurt Johansson for useful discussions, and to Alain Rouault for 
bringing our attention to the paper \cite{Pop} by I. Popescu.

A. Bufetov has been supported in part by an Alfred P. Sloan Research Fellowship, a Dynasty Foundation Fellowship, as well as an IUM-Simons
Fellowship,
by the Grant MK-6734.2012.1 of the President of the Russian Federation,
by the Programme ``Dynamical systems and mathematical control theory''
of the Presidium of the Russian Academy of Sciences,
by the RFBR-CNRS grant 10-01-93115-NTsNIL and by the RFBR grant 11-01-00654.

M. Shcherbina has been supported in part by the project "Ukrainian branch of the French-Russian Poncelet
laboratory" - "Probability problems on groups and spectral theory".

A. Soshnikov has been supported in part by the NSF grant DMS-1007558.

\section{Main results}

\subsection{Asymptotic equipartition property}
Consider a probability measure on $\mathbb{R}^N$ with the density function $P^{\beta}_N$ defined in (\ref{meas:Hermitian}).
Consider the (infinite) product of these probability spaces and denote by $\mathbb{P}^{\beta}$ the corresponding
probability measure on $\Omega:=\mathbb{R}^1\times\mathbb{R}^2\times\cdots\times\mathbb{R}^N\times\cdots$. Define random variables $X_N$ on
$\Omega$ by
\begin{equation*}
X_N(\lambda):=-\frac{\ln P^{\beta}_N(\bar{\lambda})}{N}:=-\frac{\ln P^{\beta}_N(\bar{\lambda}_N)}{N},
\end{equation*}
where $\bar{\lambda}=\{\bar{\lambda}_1,\bar{\lambda_2},\dots\}\in\Omega$.

\begin{theorem}[Asymptotic equipartition property]
\label{thm:AEP}
Let $V$ be a real analytic function  growing faster than $\log(1+\lambda^2)$, as $|\lambda|\to \infty$ whose equilibrium
density $\rho^V$ has $q$-interval support $\sigma$ $(q\ge 1)$.
Assume also that  $\rho^V$ is  generic,  which means that $\rho^V\not=0$ in the internal points
of $\sigma$, $\rho^V$ behaves like square root near the edges of $\sigma$, and the function
\begin{equation*}v(\lambda ):=2\int \log |\mu -\lambda |\rho^V (\mu )d\mu -V(\lambda )\end{equation*}
attains its maximum only if $\lambda $ belongs to  $\sigma $.
 Then
for any $\beta>0$ the random variables $X_N$ converge $\mathbb{P}^{\beta}$-almost surely to some constant $E_\beta(V).$
\end{theorem}

{\bf Remark}

{\it
In the Hermite case,}
\begin{equation}
\label{urav11}
E^{Her}_{\beta}=\ln(2\pi)-\ln\Gamma\left(1+\frac\beta2\right)+\frac{\beta}2\psi\left(1+\frac\beta 2\right)-\frac{\beta}2-\frac 12. 
\end{equation}
{\it The results of this theorem and Theorem \ref{thm:CLT} below
also hold for  the Circular, Laguerre, and Jacobi beta ensembles.  In particular,}
\begin{align*}
E^{Cir}_{\beta}&=\ln(2\pi)-\ln\Gamma\left(1+\frac\beta2\right)+\frac{\beta}2\psi\left(1+\frac\beta 2\right)-\frac{\beta}2,
\\
E^{Lag}_{\beta}&=\ln(2\pi)-\ln\Gamma\left(1+\frac\beta2\right)+\frac\beta2\psi\left(1+\frac{\beta}2\right)-\frac\beta2+\alpha\ln
\frac{\beta^2}2-\alpha-\beta\ln\beta,
\\
E^{Jac}_{\beta}&=\ln(2\pi)-\ln\Gamma\left(1+\frac\beta2\right)+\frac\beta2\psi\left(1+\frac{\beta}2\right)-\frac\beta2-(\mu+\nu)\ln 2,
\end{align*}
{\it where} $\psi$ {\it is the digamma function,}
\begin{equation}
\label{digamma}
\psi(x)=\frac{d}{dx} \log \Gamma(x).
\end{equation}
{\it For a general potential} $V,$ {\it the constant} $ E_{\beta}(V)$ {\it is defined in} (\ref{urav13}), (\ref{urav11}), {\it and} (\ref{urav12}).

{\it The fact that for the classical beta ensembles the mean energy and the specific heat
can be computed from the Selberg integral is well known in the physics community.  For the 
calculation of the specific heat for the circular beta ensemble we refer to \cite{Dys1}, \cite{Dys2}.}

The following is an immediate corollary of Theorem \ref{thm:AEP}.
\begin{corollary}
\label{cor:AEP}
For any $\beta>0$ and any $\varepsilon>0$, there exists $N_{\varepsilon,\beta}>0$ such that for any $N>N_{\varepsilon,\beta}$,
there exists $S^{\beta}_N(\varepsilon)\subset\supp(P^{\beta}_N)$ such that for any $\bar{\lambda}\in S^{\beta}_N(\varepsilon)$ we have
\begin{equation*}
e^{-(E_{\beta}+\varepsilon)N}<P^{\beta}_N(\bar{\lambda})<e^{-(E_{\beta}-\varepsilon)N}
\end{equation*}
and
\begin{equation*}
\lim_{N\rightarrow\infty} \int_{\bar{\lambda}\in S^{\beta}_N(\varepsilon)}P^{\beta}_N(\bar{\lambda}) d\bar{\lambda}=1.
\end{equation*}

\end{corollary}

\subsubsection{The Shannon-McMillan-Breiman Theorem}
The interpretation of the constant $E_\beta(V)$ as the entropy of the corresponding random matrix ensemble is
by analogy with the Shannon-McMillan-Breiman theorem. To illustrate the analogy, we briefly recall the
Shannon-McMillan-Breiman theorem in the case of the Bernoulli process.

Let $W$ be the set of sequences of $0$s and $1$s, i.e. $W=\{0,1\}^{\mathbb{N}}$. Given $\mathfrak{w}\in W$ let
$\mathfrak{w}_n=(w_1,\dots,w_n)\in\{0,1\}^n$ be the first $n$ elements of $\mathfrak{w}$, and let $S_{\mathfrak{w}_n}$
be the cyllinder set corresponding to $\mathfrak{w}_n$. $S_{\mathfrak{w}_n}$ is the set of all sequences in $W$ the
first $n$ elements of which coincide with $\mathfrak{w}_n$. Given $p\in(0,1)$, let $P_p$ be the Bernoulli measure on
$W,$ which to a cylinder set $S_{\mathfrak{w}_n}$ assigns the probability
\begin{equation*}
P_p(S_{\mathfrak{w}_n})=p^k(1-p)^{n-k},
\end{equation*}
where $k$ is the number of ones in $\mathfrak{w}_n$. Let
\begin{equation*}
H(p)=-p\ln p - (1-p)\ln(1-p)
\end{equation*}
be the entropy of the Bernoulli measure. The Shannon-McMillan-Breiman theorem states that for any $p\in(0,1)$, we have
\begin{equation*}
P_p\left(\mathfrak{w}\in W:\lim_{n\rightarrow\infty}-\frac{\ln P_p(S_{\mathfrak{w}_n})}{n}=H(p)\right)=1.
\end{equation*}
In other words, the theorem states that the random variables $-\frac{\ln P_p(S_{\mathfrak{w}_n})}{n}$ converge to the entropy $H(p)$ almost surely.

\subsection{Central Limit Theorem}

\begin{theorem}[Central Limit Theorem]
\label{thm:CLT}
Under the assumptions of Theorem \ref{thm:AEP},
$$Y_N^{\beta}:=N^{-1/2}(\log P_N^{\beta}(\bar\lambda) +N \*E_\beta(V))$$
converges in distribution as $N\to \infty$ to the Gaussian random variable with expected value $0$ and variance
\begin{equation*}
\frac{\beta}{2}-\frac{\beta^2}{4}\*\psi'\left(1+\frac \beta 2\right),
\end{equation*}
where $\psi$ is defined in (\ref{digamma}).
\end{theorem}

{\bf Remark} 

{\it The results of Theorems \ref{thm:AEP} and \ref{thm:CLT} in the special case of the Hermite 
beta ensemble (\ref{meas:beta}) were first proven by Ionel Popescu in \cite{Pop}.}

{\bf Remark}

{\it Let us write the probability density (\ref{meas:Hermitian}) as}
\begin{equation*}
P^{\beta}_N(\bar{\lambda})=\frac{e^{-\beta H_N(\bar{\lambda})/2}}{Z_N(\beta)},
\end{equation*}
{\it where the potential energy term} $H_N$ {\it is given by}
\begin{equation}
\label{potenergy}
H_N(\bar{\lambda})=N\*\sum_{i=1}^N V(\lambda_i) - \sum_{i \not=j} \ln |\lambda_i-\lambda_j|.
\end{equation}
{\it Then the Central Limit Theorem result can be reformulated in terms of convergence in distribution of}
\begin{equation*}
W_N(\bar{\lambda})=\frac{H_N(\bar{\lambda})-C_{N,\beta}}{\sqrt{N}}
\end{equation*}
{\it to the Gaussian limit} $N(0,\frac 2{\beta}-\psi'\left(1+\frac \beta 2\right))$, {\it where}
$C_{N,\beta}$ { \it is an appropriate centering constant.

In particular, for the classical beta ensembles, one has}
\begin{align*}
C^{Her}_{N,\beta}&=\frac 38 N^2-\frac 12 N\ln N+\left(-\frac 12 \ln\frac \beta 2-\frac 14+\frac 12\psi\left(1+\frac \beta 2\right)\right)N,
\\
C^{Cir}_{N,\beta}&=-\frac 12 N\ln N-\left(\frac 12 \ln\frac \beta 2-\frac 12\psi\left(1+\frac \beta 2\right)\right)N,
\\
C^{Lag}_{N,\beta}&=\left(\frac 34+\frac{\ln 2}{2}\right)N^2-\frac 12 N\ln N-\frac{1+\ln\beta-\psi\left(1+\frac{\beta}2\right)}2 N,
\\
C^{Jac}_{N,\beta}&=\frac{\ln 2}2 (N-2)N-\frac 12 N\ln N+\left(-\frac 12 \ln\frac \beta 2+\frac 12\psi\left(1+\frac \beta 2\right)\right)N.
\end{align*}

{\bf Remark}

{\it It should be noted that both terms in the expression (\ref{potenergy}) for the potential energy, namely}
$N\*\sum_{i=1}^N V(\lambda_i)$ and $\sum_{1\leq i< j \leq N} \ln |\lambda_i-\lambda_j|,$
{\it have 
fluctuations of order} $N,$ { \it (see e.g. \cite{JohDuke1998}, \cite{Pastur}).  
At the same time, the sum of these two terms fluctuates on a much smaller order, namely}
$N^{1/2}.$ {\it The cancellations take place
because of the Euler-Lagrange variational equations (\ref{EL1}-\ref{EL2}) for the equilibrium measure. As a result, the difference}
$V(\lambda_i) - 2\*\sum_{j:j\not=i} \log|\lambda_j-\lambda_i| $ {\it is equal (up to a negligible error term) to a constant independent of}
$ \ 1\leq i\leq N.$

{\it It is also important to note that while the fluctuation of the linear statistic} 
$\sum_{i=1}^N V(\lambda_i)$ {\it is (asymptotically) Gaussian in the one-cut case, it is 
non-Gaussian, in general, in the multi-cut case (see \cite{Pastur}).}

For the Hermite beta ensemble, the potential energy $H_N(\bar{\lambda})$ attains its (unique) minimum at the configuration given by the
zeroes of the (rescaled) Hermite polynomial of order $N$, namely $h_N(\sqrt{\frac{N}{2}}\*x),$ where
$$h_N(x)=N! \* \sum_{m=0}^{[N/2]} \frac{(-1)^m \*(2x)^{N-2m}}{m!\*(N-2m)!}$$
(see e.g. \cite[A.6]{mehta2004random}).
The minimum value of $H_N$ is equal to
\begin{equation}
\label{dig}
\frac{N \*(N-1)}{4}\*(\log N +1) - \frac{1}{2} \*\sum_{j=1}^N j \* \log j=\frac{3}{8}\*N^2 -\frac{1}{2}\*N\*\log N
-\frac{1}{4}\*N +o(N).
\end{equation}
Note that the inequality $\psi(1+x) \geq \log x$ implies that
$$C^{Her}_{N,\beta} \geq \frac{3}{8}\*N^2 -\frac{1}{2}\*N\*\log N -\frac{1}{4}\*N.$$
Thus, for a typical configuration the difference between the potential energy and its minimal value is proportional to $N.$

Similar explicit computations could be done for the other classical beta ensembles.
In particular, for the Circular beta ensemble, the ground state is unique up to a rotation and is given by the vertices of a regular $N-$gon,
and the minimum of the potential energy is equal to $-\frac{1}{2}\*N \*\log N.$  For Laguerre and Jacobi beta ensembles, the ground state is unique and
is given by the zeroes of the corresponding (rescaled) orthogonal polynomial of degree $N.$   Since the computations are very similar to the ones in
\cite[A.6]{mehta2004random}), we leave the details to the reader.

\section{Proof of Theorem \ref{thm:CLT}}
We first give the proof in the case of the Hermite beta ensemble.  Then we quickly extend it to analytic potentials $V$ in
both  one-cut and multi-cut cases.

We start with the Selberg type integral corresponding to the Hermite orthogonal polynomials \cite[17.6.7]{mehta2004random}
\begin{align}
\label{eq:PartF} Z^{Her}_N(\beta)=(2\pi)^{\frac N2}
\left(\frac{N\beta}2\right)^{-\frac N2((N-1)\beta/2+1)}\prod_{j=1}^N\frac{\Gamma(1+j\beta/2)}{\Gamma(1+\beta/2)},
\end{align}
which can be rewritten as
\begin{align*}
\log Z_N^{(b)}(\beta)=&\frac{N}{2}\log 2\pi-\Big(\frac{\beta N^2}{4}+\frac{N}{2}\big(1-\frac{\beta}{2}\big)\Big)
\log\frac{\beta N}{2}-N\log \Gamma(1+\beta /2)\\&+
\sum_{j=1}^N\log \Gamma(1+\beta j/2).
\end{align*}
Let us
represent $\sum_{j=1}^N\log \Gamma(1+\beta j/2)$ in terms of the Barnes G-function defined in
\cite{Barnes} (see also \cite{Fo:10} formula (4.181)). It  satisfies the functional equation
$$ G(1+z)=\Gamma(z)G(z),\quad G(1)=1,$$
and so at the integer points $G$ can be represented as
$ G(1+N)=\prod_{j=1}^{N-1}\Gamma(1+j).$
The asymptotic expansion for $\log G$ is known (see \cite{Fo:10} formula (4.184)):
\begin{equation}
\label{Barn}
\log G(1+z)=\frac{z^2}{2}\log z-\frac{3}{4}z^2+\frac{z}{2}\log 2\pi-\frac{1}{12}\log z+\zeta'(-1)+o(1),
\quad z\to\infty.
\end{equation}
To obtain the representation for $\sum_{j=1}^N\log \Gamma(1+\beta j/2)$ we use the Stirling formula in the form
\begin{equation}\label{St}
\log \Gamma(1+z)=\big(z+\frac{1}{2}\big)\log z-z+
\frac{1}{2}\log 2\pi+\frac{1}{12 z}+r(z)\quad|r(z)|\le Cz^{-2},\end{equation}
where the bound for $r(z)$ is uniform for any interval $(\delta,\infty)$ with $\delta>0$. Then
\begin{align}
\sum_{j=1}^N\log \Gamma&(1+\beta j/2)=\sum_{j=1}^N\Big(\big(\frac{\beta j}{2}+\frac{1}{2}\big)\log\frac{\beta j}{2}-\frac{\beta j}{2}+
\frac{1}{2}\log 2\pi+\frac{1}{6\beta j}+O((\beta j)^{-2})\Big)\notag\\=&
\sum_{j=1}^N\Bigg(\Big(\frac{\beta}{2}\big(j+\frac{1}{2}\big)\log j-\frac{\beta j}{2}+
\frac{\beta}{4}\log 2\pi+\frac{\beta}{12 j}\Big)+r(\beta j) \notag\\&+
\frac{1}{2}\big(1-\frac{\beta}{2}\big)\log j+\big(\frac{\beta j}{2}+\frac{1}{2}\big)\log\frac{\beta }{2}+
\frac{1}{2}\big(1-\frac{\beta}{2}\big)\log 2\pi-\frac{1}{12 j}\big(\frac{\beta}{2}-\frac{2}{\beta}\big)\Bigg)\notag\\
=&\sum_{j=1}^N\frac{\beta}{2}\log\Gamma(1+ j)+\frac{1}{2}\big(1-\frac{\beta}{2}\big)\log\Gamma(1+ N)+
\big(\frac{\beta N^2}{4}+\frac{N}{2}\big(1+\frac{\beta}{2}\big)\big)\log\frac{\beta}{2}\notag\\&+
\frac{N}{2}\big(1-\frac{\beta}{2}\big)\log 2\pi-\frac{1}{12}\big(\frac{\beta}{2}-\frac{2}{\beta}\big)\log N+C_N(\beta)
+o(1)\notag\\=&
\frac{\beta}{2}\log G(1+N)+\frac{1}{2}\big(1+\frac{\beta}{2}\big)\log\Gamma(1+ N)
+\big(\frac{\beta N^2}{4}+\frac{N}{2}\big(1+\frac{\beta}{2}\big)\big)\log\frac{\beta}{2}\notag\\&+
\frac{N}{2}\big(1-\frac{\beta}{2}\big)\log 2\pi-\frac{1}{12}\big(\frac{\beta}{2}-\frac{2}{\beta}\big)\log N+C(\beta)+o(1).
\label{pr_Gamma}\end{align}
Here
\begin{align*}
& C(\beta)=\lim_{N\to\infty}C_N(\beta), \\
& C_N(\beta)=\sum_{j=1}^N\Big(r(\beta j/2)-\frac{\beta}{2}r(j)\Big)
-\frac{1}{12}\big(\frac{\beta}{2}-\frac{2}{\beta}\big)\gamma,
\end{align*}
where $\gamma$ denotes the Euler constant and  $r(z)$ is defined in (\ref{St}).

This expansion combined with (\ref{Barn}) and (\ref{St}) implies
\begin{align*}
\sum_{j=1}^N\log \Gamma(1+\beta j/2)=&\frac{\beta N^2}{4}\log\frac{\beta N}{2}-\frac{3}{8}\beta N^2+
\frac{N}{2}\big(1+\frac{\beta}{2}\big)\log \frac{\beta N}{2}-\frac{N}{2}\big(1+\frac{\beta}{2}\big)
\\&+
\frac{N}{2}\log 2\pi+R(\beta)\log N+\widetilde C(\beta)+o(1),
\end{align*}
where
\begin{align*}
& R(\beta)=\frac{\beta}{24}+\frac{1}{4}+\frac{1}{6\beta}, \\
& \widetilde C(\beta)=C(\beta)
+\frac{1}{4}(1+\frac{\beta }{2})\log 2\pi+\zeta'(-1).
\end{align*}

Applying the Selberg formula (\ref{eq:PartF}),
the above representation yields
\begin{align}\notag
\log Z_N^{Her}(\beta)=&\frac{\beta}{2} N^2\mathcal{E}_0+\frac{\beta N}{2}\log\frac{\beta N}{2}+N\Big(\log 2\pi-\log \Gamma(1+\frac{\beta }{2})
-\frac{1}{2}(1+\frac{\beta }{2})\Big)\\&+
R(\beta)\log N+\widetilde C(\beta)+o(1),\quad \mathcal{E}_0=-\frac{3}{4}.
\label{logZ_0}
\end{align}

Consider now the random variable
$$ Y_N^{Her,\beta}=N^{-1/2}(\log P_N^{Her,\beta}(\bar\lambda)+N \*E_\beta^{Her})
=N^{-1/2}\Big(-\frac{\beta}{2} H_{N}(\bar\lambda)-\log Z_N(\beta)+N\*E_\beta^{Her}\Big).$$
The Laplace transform of the distribution of $Y_N^{Her,\beta},$
\begin{equation}
\Phi(t)=\int d\bar\lambda  e^{-t\*Y_N^{Her,\beta}}\* P_N^{Her,\beta}(\bar\lambda) ,
\end{equation}
can be written as
\begin{equation}
\Phi(t)=( Z_N^{Her}(\beta))^{-1}\int d\bar\lambda\exp\Big\{-\frac{\beta}{2}(1+\frac{t}{\sqrt{N}})H_N(\bar\lambda)
-\frac{t}{\sqrt{N}}\log Z_N(\beta)+ t\sqrt{N}E_\beta^{Her}\Big\}.
\end{equation}
Thus,
\begin{align*}
\log \Phi(t)=&\log Z_N(\beta(1+\frac{t}{\sqrt{N}}))-(1+\frac{t}{\sqrt{N}})\log Z_N(\beta)+\beta t\sqrt{N}E_\beta^{Her}.
\end{align*}
Since according to (\ref{logZ_0})
\begin{align}
\label{logZ_g}
\log Z_N^{Her}(\beta)=\frac{\beta}{2} N^2\mathcal{E}_0+\frac{\beta N}{2}\log N
+Nf_0(\beta)+R(\beta)\log N+\widetilde C(\beta)+o(1)
\end{align}
with
$$ f_0(\beta)=\log 2\pi-\log \Gamma(1+\frac{\beta }{2})+\frac{\beta }{2}\log\frac{\beta }{2}
-\frac{1}{2}(1+\frac{\beta }{2}),$$
one can see immediately that
$$\log Z_N(\beta(1+\frac{t}{\sqrt{N}}))-(1+\frac{t}{\sqrt{N}})\log Z_N(\beta)=t\sqrt N(\beta f'_0(\beta)-f_0(\beta))+
\frac{\beta^2t^2}{2}f''(\beta)+o(1).$$
Hence, using  that by definition of $E_\beta^{Her}$ in (\ref{urav11})
$$ E_\beta^{Her}=f_0(\beta)-\beta f'_0(\beta),$$
we get
$$\Phi(t)=\exp\{\frac{\beta^2t^2}{2}f''(\beta)+o(1)\}.$$
Thus, we have proved that $Y_N^{Her,\beta}$ converges in distribution to the Gaussian random
variable with zero mean and variance
$$\beta^2 \*f''(\beta)=\frac{\beta}{2}-\frac{\beta^2}{4}\psi'(1+\frac{\beta}{2}).$$

Consider now the case of general one-cut potential $V$.  We write
\begin{align}
\label{H}
H_N(\bar\lambda;V)=&-N\sum_{i=1}^n
V(\lambda_i)+\sum_{i\not=j}\log|\lambda_i-\lambda_j|,\\
P_N^{\beta}(\bar\lambda;V)=&(Z_N(\beta;V))^{-1}e^{-\beta H_N(\bar\lambda;V)/2},\\
Z_N^{(b)}(\beta;V)=&\int e^{-\beta H_N(\bar\lambda;V)/2}d\bar\lambda.
\end{align}
 We use the expansion for $\log (Z_N(\beta;V)/N!)$ found in \cite{S:11} (formulas
(1.31)-(1.33)) in the modified form of \cite{S:12} (formulas (1.20), (1.21)). Note that $F_\beta(N)$ below differs from
$F_\beta$ of (1.21)  by $\log N!$.  We have
\begin{align}\label{logZ[V]}
\log Z_N(\beta;V)=&\frac{\beta N^2}{2}\mathcal{E}[ V]+F_\beta(N)
+NS[\rho^V]\Big(\frac{\beta}{2}-1\Big)+r_\beta[ \rho^V]+O(N^{-1}),\\
F_\beta(N)=&\log Z_N^{Her}(\beta)-\frac{\beta N^2}{2}\mathcal{E}_0,\notag\\
\label{urav12}
S[\rho^V]=&\big((\log \rho^V,\rho^V)-1+\log 2\pi\big),\\
\mathcal{E}[ V]=&-(V,\rho^V)+(L\rho^V,\rho^V),\quad (Lf)(\lambda)=\int \log |\lambda-\mu|f(\mu) d\mu,\notag
\end{align}
where $\rho^V$ is  the equilibrium density,
$r_\beta[\rho]$ is a smooth function of $\beta$ defined by a contour integral with the Stieltjes transform
of $\rho^V,$ and $(f,g):=\int f(x)\*g(x)\*dx.$
Since (\ref{logZ[V]}) is also written in the form (\ref{logZ_g}) with $\mathcal{E}_0$ replaced by
$\mathcal{E}[ V]$, and $f_0$, $\widetilde C(\beta)$
replaced by
$$ f_V(\beta)=f_0(\beta)+S[\rho^V]\Big(\frac{\beta}{2}-1\Big),\quad C[V]=\widetilde C(\beta)+r_\beta[ \rho^V],$$
by the same argument we conclude that the CLT is valid for
$$Y_N^{\beta}=N^{-1/2}(\log P_N^{\beta}(\bar\lambda;V)+N \*E_\beta[V])
=N^{-1/2}\Big(-\frac{\beta}{2} H_{N}(\bar\lambda;V)-\log Z_N(\beta;V)+N \*E_\beta[V]\Big),$$
where
\begin{equation}
\label{urav13}
E_\beta[V]=E_\beta^{Her}-S[\rho^V],
\end{equation}
and the variance is the same as in the Gaussian case. 

We are left to consider the multi-cut case. Since the analog
of (\ref{logZ[V]}) obtained in \cite{S:12} in the multi-cut case is more complicated, we first need some extra definitions.

Denote
\begin{equation}\label{mu^*}
  \sigma:=\bigcup_{\alpha=1}^q\sigma_\alpha,\quad \sigma_\alpha=[a_{\alpha},b_{\alpha}],\quad
    \mu_\alpha:=\int_{\sigma_\alpha}\rho_\alpha^V(\lambda)d\lambda,\quad \rho_\alpha^V:=\mathbf{1}_{\sigma_\alpha}\rho^V.
\end{equation}

  Define the operators $\mathcal{L}$, $\widehat{\mathcal{L}}$ and $\widetilde{\mathcal{L}}$ by
\begin{equation}\label{cal-L}
    \mathcal{L}f:=\mathbf{1}_\sigma L[f\mathbf{1}_\sigma],\quad\widehat{\mathcal{L}}_\alpha f:=
    \mathbf{1}_{\sigma_\alpha} L[f\mathbf{1}_{\sigma_\alpha}],\quad \widehat{\mathcal{L}}:=
   \oplus_{\alpha=1}^q\widehat{\mathcal{L}}_\alpha,\quad
    \widetilde{\mathcal{L}}:=\mathcal{L}-\widehat{\mathcal{L}},
\end{equation}
on the set of the functions
\begin{equation*} \mathcal{H}=\oplus_{\alpha=1}^q L_1[\sigma_\alpha].\end{equation*}
Note that the topology of $\mathcal{H}$ is not important  below.

For  each interval $\sigma_\alpha$ we also define the operator
\begin{align}&\overline D_\alpha=\frac{1}{2}(D_\alpha+D^*_\alpha),\quad D_\alpha h(\lambda)=
\frac{1}{\pi^2 }\int_{\sigma_\alpha}
\frac{h'(\mu)\sqrt{(\mu-a_\alpha)(b_\alpha-\mu)} d\mu}{(\lambda-\mu)\sqrt{(\lambda-a_\alpha)(b_\alpha-\lambda)}},\label{bar_D}\end{align}
and set
\begin{equation}\label{D,L}
   \overline D:=\oplus_{\alpha=1}^q \overline D_\alpha.
\end{equation}

Let $\mathcal{Q}$ be
a positive definite  $q\times q$ matrix of the form
\begin{equation}\label{Q}
\mathcal{Q}=\{\mathcal{Q}_{\alpha\alpha'}\}_{\alpha,\alpha'=1}^{q},\quad
\mathcal{Q}_{\alpha\alpha'}=-(\mathcal{L}\psi^{(\alpha)},\psi^{(\alpha')}),\end{equation}
where the function $\psi^{(\alpha)}(\lambda)$ is a unique solution of the system of equations
\begin{equation}\label{cond_psi}
-(\mathcal{L}\psi^{(\alpha)})_{\alpha'}=\delta_{\alpha\alpha'},\quad \alpha'=1,\dots,q.
\end{equation}
It is known that $\psi^{(\alpha)}$ can be chosen in the form
$$\psi^{(\alpha)}(\lambda)=
p_\alpha(\lambda) X^{-1/2}_\sigma(\lambda)\mathbf{1}_\sigma,\quad X^{1/2}_\sigma(\lambda)=
\Im\Big(\prod (z-a_\alpha)(z-b_\alpha)\Big)^{1/2}\Big|_{z=\lambda+i0},$$
where $p_\alpha$ is a polynomial of degree $q-1$.
Set
\begin{equation}\label{I[h]}I[h]=(I_1[h],\dots,I_q[h]),\quad I_\alpha[h]:=
\sum_{\alpha'}\mathcal{Q}^{-1}_{\alpha\alpha'}(h,\psi^{(\alpha')}),\end{equation}
and define a quasi-periodic in $N$ function $\Theta(\beta;\{N\bar\mu\})$ by
\begin{align}\label{theta}
&\Theta( \beta;\{N\bar\mu\}):=\sum_{n_1+\dots+n_q=n_0} \exp\Big\{-\frac{\beta}{2}\Big(\mathcal{Q}^{-1}\Delta\bar n,\Delta\bar n\Big)
+\big(\frac{\beta}{2}-1\big)( \Delta\bar n, I[\log\overline\rho^V])\Big\},\\
&\{N\bar\mu\}=(\{N\mu_1\},\dots,\{N\mu_q\}),\quad(\Delta\bar  n)_\alpha=n_\alpha-\{N\mu_\alpha\},
\quad n_0=\sum_{\alpha=1}^q\{N\mu_\alpha\},
\notag\end{align}
with  $\mathcal{Q}$ of (\ref{Q}), $\{\mu_\alpha\}_{\alpha=1}^q$ of (\ref{mu^*}), $I[h]$ of (\ref{I[h]}), and
 $\log\overline\rho=(\log\rho_1,\dots,\log\rho_q)$.

Using the above definitions we can write now the analog of (\ref{logZ[V]}) obtained in \cite{S:12}, formula (1.34):
\begin{align}\label{as_Q}
\log Z_N^{(b)}[\beta,V]=&
\frac{N^2\beta}{2}\mathcal{E}[V]+F_\beta(N)+N(\frac{\beta}{2}-1)S[\rho^V]+(q-1)\big(R(\beta)\log N+\widetilde C(\beta)\big)
\notag\\&
+\sum_{\alpha=1}^q(r_{\beta}[\mu_\alpha^{-1}\rho_\alpha]+R(\beta)\log\mu_\alpha)-
\frac{1}{2}\log\det(1-\overline D\widetilde{\mathcal{L}})\\
&+\dfrac{2}{\beta}\Big(\dfrac{\beta}{2}-1\Big)^2
\big(\widetilde{\mathcal{L}}(1-\overline D\widetilde{\mathcal{L}})^{-1}\overline\nu,\overline\nu\big)
+\log\Theta(\beta;\{N\bar\mu\})+O(N^{-\kappa}), \quad\kappa>0,
\notag\end{align}
where $\mu_\alpha,\,\rho_\alpha$ are defined in (\ref{mu^*}),  $r_{\beta}[\mu_\alpha^{-1}\rho_\alpha]$ for each
$\sigma_\alpha$ is the same as in (\ref{logZ[V]}), $F_\beta(N)$, $R(\beta)$, and $\widetilde C(\beta)$ are defined in
(\ref{logZ_0}) (note that $R(\beta)$ corresponds to $-c_\beta$ from (1.34) of \cite{S:12}, $\widetilde C(\beta)$
corresponds to $c_\beta^{(1)}$ and $F_\beta(N)$ differs from that of (1.34) by $\log N!$), and $\det$ means the
Fredholm determinant of $\overline D\widetilde{\mathcal{L}}$ on $\sigma$. The
non positive measures $(\nu_1,\dots,\nu_q)=:\overline\nu$ have the form
\begin{align*}
(\nu_\alpha,h):=&\frac{1}{4}(h(b_\alpha)+h(a_\alpha))-\frac{1}{2\pi}\int_\sigma\frac{h(\lambda)d\lambda}
{\sqrt{(\lambda-a_\alpha)(b_\alpha-\lambda)}}
+\frac{1}{2}(D_\alpha\log P_\alpha,h).
\end{align*}
with
\begin{equation*}P_\alpha(\lambda)=\frac{2\pi\rho^V_\alpha(\lambda)}{\sqrt{(\lambda-a_\alpha)(b_\alpha-\lambda)}}.\end{equation*}
It is easy to see that in the multi-cut case the expression (\ref{as_Q}) still can be written in the form (\ref{logZ_g}) with
$\mathcal{E}_0$ replaced by
$\mathcal{E}[ V]$, $R(\beta)$ replaced by $qR(\beta)$, and $f_0$, $\widetilde C(\beta)$
replaced by
\begin{align*} f_V(\beta)=&f_0(\beta)+S[\rho^V]\Big(\frac{\beta}{2}-1\Big),\\
C_q[V]=&q\widetilde C(\beta)+\sum_{\alpha=1}^q(r_{\beta}[\mu_\alpha^{-1}\rho_\alpha]+R(\beta)\log\mu_\alpha)-
\frac{1}{2}\log\det(1-\overline D\widetilde{\mathcal{L}})\\
&+\dfrac{2}{\beta}\Big(\dfrac{\beta}{2}-1\Big)^2
\big(\widetilde{\mathcal{L}}(1-\overline D\widetilde{\mathcal{L}})^{-1}\overline\nu,\overline\nu\big)
+\log\Theta(\beta;\{N\bar\mu\}).\end{align*}
Since the structure of $\Theta(\beta;\{N\bar\mu\})$ (see (\ref{theta})) guarantees that
$$ \Big|\log\Theta\Big(\beta;\{N\bar\mu\}\Big)-\log\Theta\Big(\beta\Big(1+\frac{t}{\sqrt N}\Big);\{N\bar\mu\}\Big)\Big|
\le \frac{K(\beta) |t|}{\sqrt N},$$
by the same argument we conclude that the CLT is valid for
$$Y_N^{\beta}=N^{-1/2}(\log P_N^{\beta}(\bar\lambda;V)+N \*E_\beta[V])
=N^{-1/2}\Big(\frac{\beta}{2} H_{N}(\bar\lambda;V)-\log Z_N(\beta;V)+N \*E_\beta[V]\Big),$$
where
$$E_\beta[V]=E_\beta^{Her}-S[\rho^V],$$
and the variance is the same as in the Gaussian case.

The theorem is proven.

{\bf Remark}

{\it As we have mentioned above, the result of Theorem \ref{thm:CLT} also holds for Circular, Laguerre, and Jacobi beta ensembles.
The proofs are very similar to the one in the Hermite case.  The arguments are
based on the explicit formulas for the partition functions given by the Selberg integrals}
\cite{mehta2004random}:
\begin{align*}
Z^{Cir}_N(\beta)=&(2\pi)^{N}\frac{\Gamma(1+N\beta/2)}{\Gamma(1+\beta/2)^N},\\
Z^{Lag}_N(\beta)=&(\beta N)^{-\frac{N(N-1)}{2}\beta-\alpha N} \prod_{j=0}^{N-1}\frac{\Gamma(1+(j+1)\beta/2)\Gamma(\alpha+j\beta/2)}
{\Gamma(1+\beta/2)}, \\
Z^{Jac}_N(\beta)=&2^{\frac \beta 2 N(N-1)+N(\mu+\nu-1)}
\prod_{j=0}^{N-1}\frac{\Gamma(1+(j+1)\beta/2)\Gamma(\mu+j\beta/2)\Gamma(\nu+j\beta/2)}{\Gamma(1+\beta/2)\Gamma(\mu+\nu+(N+j-1)\beta/2)}.
\end{align*}
{\it The details are left to the reader.}

To prove Theorem \ref{thm:AEP} it suffices to observe that by the Chebyshev inequality we have
$$ Prob\{|X_N^{\beta}+E_\beta|>\varepsilon\}=Prob\{|Y_N^{\beta}|>\sqrt N\varepsilon\}\le
(\Phi(1)+\Phi(-1))e^{-\sqrt N\varepsilon}\le Ce^{-\sqrt N\varepsilon},$$
and apply the Borel-Cantelli lemma.

\end{document}